%% file: paper.tex
\documentclass[10pt, conference]{IEEEtran}

\newcommand{\ignore}[1]{#1}
\newcommand{\algfontsize}{\small}

\ifCLASSOPTIONcompsoc
  \usepackage[nocompress]{cite}
\else
  \usepackage{cite}
\fi
\ifCLASSINFOpdf
\else
\fi
\ifCLASSOPTIONcompsoc
  \usepackage[caption=false,font=footnotesize,labelfont=sf,textfont=sf]{subfig}
\else
  \usepackage[caption=false,font=footnotesize]{subfig}
\fi
\usepackage{url}

\usepackage{amsfonts,amsmath,latexsym,amssymb}
\usepackage{algorithm,algpseudocode}
\usepackage{graphicx,xcolor}
\usepackage{hyperref,cleveref}
\usepackage{tikz}
\usetikzlibrary{3d}
\usepackage{pgfplots}
\usepackage{pgfplotstable}
\usepackage{etoolbox}
\usetikzlibrary{patterns}
\usepackage{ourmacros}

\begin{document}
%
\title{Parallel Nonnegative CP Decomposition of~Dense~Tensors}

\author{\IEEEauthorblockN{Grey Ballard and Koby Hayashi}
\IEEEauthorblockA{Wake Forest University\\
Winston Salem NC 27109 \\
Email: \{ballard,hayashi\}@wfu.edu}
\and
\IEEEauthorblockN{Ramakrishnan Kannan}
\IEEEauthorblockA{Oak Ridge National Laboratory\\
Oak Ridge, TN 37830\\
Email: kannanr@ornl.gov}}



%


\maketitle

\input{abstract}


%
\IEEEpeerreviewmaketitle

\input{introduction}

\input{prelims}

\input{survey}

\input{algorithm}

\input{experiments}

\input{conclusion}


\ifCLASSOPTIONcompsoc
\else
\fi



%


\ignore{
\newpage 
\input{appendix}

}


\end{document}

%% file: abstract.tex

\begin{abstract}

The CP tensor decomposition is a low-rank approximation of a tensor. 
We present a distributed-memory parallel algorithm and implementation of an alternating optimization method for computing a CP decomposition of dense tensor data that can enforce nonnegativity of the computed low-rank factors.
The principal task is to parallelize the matricized-tensor times Khatri-Rao product (MTTKRP) bottleneck subcomputation.
The algorithm is computation efficient, using dimension trees to avoid redundant computation across MTTKRPs within the alternating method.
Our approach is also communication efficient, using a data distribution and parallel algorithm across a multidimensional processor grid that can be tuned to minimize communication.
We benchmark our software on synthetic as well as hyperspectral image and neuroscience dynamic functional connectivity data, demonstrating that our algorithm scales well to 100s of nodes (up to 4096 cores) and is faster and more general than the currently available parallel software. 

\end{abstract}

%% file: introduction.tex

\section{Introduction}

The CP decomposition is a low-rank approximation of a multi-dimensional array, or tensor, which generalizes matrix approximations like the truncated singular value decomposition.
It approximates the input tensor by a sum of rank-one tensors, which are outer products of vectors.
CP is often used for finding hidden patterns, or latent factors, within tensor data, particularly when the goal is to interpret the factors, and it is popular within the signal processing, machine learning, and scientific computing communities.


To aid in interpretability, domain-specific constraints are often imposed on the computed factors.
We focus in this paper on dense tensors (when nearly all of the tensor entries are nonzero) and on constraining solutions to have nonnegative entries, which is useful when the tensor data itself is nonnegative. Formally, NNCP can be defined as
\SplitN{\label{eqn:nncp}}{
 \min_{\HH^{(i)} \geq 0} & \left\| \TA - \sum_{r=1}^R \Mn{H}{1}(:,r) \circ \cdots \circ \Mn{H}{N}(:,r) \right\|^2
}
where $\Mn{H}{1}(:,i) \circ \cdots \circ \Mn{H}{N}(:,i)$ is the outer product of the $i^{th}$ vector from all the $N$ factors that yields a rank one tensor $\T{M}$ and $\sum_{r=1}^R \Mn{H}{1}(:,r) \circ \cdots \circ \Mn{H}{N}(:,r)$ results in a sum of $R$ rank one tensors that will be of the same dimension as the input tensor $\TA$. 
For example, in imaging and microscopy applications, tensor values often correspond to intensities, and NNCP can be used to cluster and analyze the data in a lower-dimensional space \cite{JC+16}.
In this work, we consider two such applications: a series of time-lapse hyperspectral images \cite{FAN16} and a dynamic functional correlation data set generated from functional magnetic resonance images of human brains \cite{VEU+12}.

One approach to handling multidimensional data is to ``matricize'' it, combining sets of modes to reshape the data into a matrix, so that standard matrix  methods like principal component analysis or nonnegative matrix factorization can be applied.
While this approach can be effective in certain cases, reshaping the data destroys multidimensional relationships among entries that the matrix methods cannot recover.
By maintaining the tensor structure of the data, the low-rank representations preserve these relationships, often producing better and more interpretable results.

However, tensor methods are more complicated both mathematically and computationally.
The kernel computations within standard algorithms for computing NNCP can be formulated as matrix computations, but the complicated layout of tensors in memory prevents the straightforward use of BLAS and LAPACK libraries.
In particular, the matrix formulation of subcomputations involve different views of the tensor data, so no single layout yields a column- or row-major matrix layout for all subcomputations.
Likewise, the parallelization approach for tensor methods is not a straightforward application of parallel matrix computation algorithms.

In developing an efficient parallel algorithm for computing a NNCP of a dense tensor, the key is to parallelize the bottleneck computation known as Matricized-Tensor Times Khatri-Rao Product (MTTKRP), which is performed repeatedly for each mode of the tensor.
The parallelization must load balance the computation, minimize communication across processors, and distribute the results so that the rest of the computation can be performed independently.
In our algorithm, not only do we load balance the computation, but we also compute and store temporary values that can be used across MTTKRPs of different modes using a technique known as dimension trees, significantly reducing the computational cost compared to standard approaches.
Our parallelization strategy also avoids communicating tensor entries and minimizes the communication of factor matrix entries, helping the algorithm to remain computation bound and scalable to high core counts.

As we detail in the related work, the general techniques for reducing computation and communication have been used in similar contexts.
The recomputation avoidance was proposed in a sequential algorithm \cite{PTC13a}, the parallelization scheme was proposed and analyzed for general tensors \cite{BKR17-TR}, and the algorithm was implemented for 3D tensors \cite{LK+17b}.

We summarize our main contributions as follows:
\begin{itemize}
		\item we present the first distributed-memory parallel implementation of NNCP algorithms for arbitrary-dimension dense tensors,
		\item we optimize the use of dimension trees for dense tensors, avoiding recomputation across multiple MTTKRPs,
		\item our parallel algorithm is communication optimal with a carefully chosen processor grid,
		\item we demonstrate a performance improvement of up to $2.2\times$ over the existing state-of-the-art parallel software on 3D tensors,
		\item our implementation obtains efficient parallel scaling of up to $1771\times$ on 4096 cores.
\end{itemize}

%% file: prelims.tex

\section{Preliminaries} 
\label{sec:prelims}

\subsection{Notation}

Tensors will be denoted using Euler script (e.g., $\T{T}$), 
matrices will be denoted with uppercase boldface (e.g., $\M{M}$), vectors will be denoted with lowercase boldface  
(e.g., $\V{v}$), and scalars will not be boldface (e.g., $s$). 
We use Matlab style notation to index into tensors, matrices, and vectors, and we use 1-indexing. 
For example, $\M{M}(:,c)$ gives the $c$th column of the matrix M.

We use $\circ$ to denote the outer product of two or more vectors.
The Hadamard product is the element-wise matrix product and will be denoted using $\Hada$. 
The Khatri-Rao product, abbreviated KRP, will be denoted with $\Khat$. 
Given matrices $\M{A}$ and $\M{B}$ that are $I_{A} \times R$ and $I_{B} \times R$, the KRP $\M{K} = \M{A} \Khat \M{B}$ is $I_{A}I_{B} \times R$. 
It can be thought of as a row-wise Hadamard product, where $\M{K}(i+I_{A}(j{-}1),:) = \M{A}(i,:) \Hada \M{B}(j,:)$, or a column-wise Kronecker product, where $\M{K}(:,c) =  \M{A}(:,c) \Kron \M{B}(:,c)$.

The CP decomposition of a tensor (also referred to as the CANDECOMP/PARAFAC or canonical polyadic decomposition) is a low-rank approximation of a tensor, where the approximation is a sum of rank-one tensors and each rank-one tensor is the outer product of vectors.
We use the notation
$$\TA \approx \CP = \sum_{r=1}^R \Mn{H}{1}(:,r) \circ \cdots \circ \Mn{H}{N}(:,r)$$
to represent a rank-$R$ CP model, where $\Mn{H}{n}$ is called a factor matrix and collects the mode-$n$ vectors of the rank-one tensors as columns.
The columns of the factor matrices are often normalized, with weights collected into an auxiliary vector $\V{\lambda}$ of length $R$; in this case we use the notation $\CPl$.

A Nonnegative CP decomposition (NNCP) constrains the factor matrices to have nonnegative values.
In this work, we are interested in NNCP models that are good approximations to $\TA$ in the least squares sense.
That is, we seek 
$$\min_{\HH^{(i)} \geq 0}  \| \TA - \CPl \|,$$
where the tensor norm is a generalization of the matrix Frobenius or vector 2-norm, the square root of the sum of squares of the entries.

The $n$th mode matricized tensor denoted by $\M{A}_{(n)}$ is a $I_n\times I/I_n$ matrix formed by organizing the $n$th mode fibers of a tensor $\T{X}$ with dimensions $I_{1} \times ... \times I_{N}$ (and $I=\prod I_n$) into the columns of a matrix. 
The Matricized-Tensor Times Khatri-Rao Product or MTTKRP will be central to this work and takes the form $\Mn{M}{n} = \Mz{A}{n} \Mn{K}{n}$, where $\Mn{K}{n}= \Mn{H}{N} \Khat \cdots \Khat \Mn{H}{n+1} \Khat \Mn{H}{n-1} \Khat  \cdots \Khat \Mn{H}{1}$. 

\subsection {Block Coordinate Descent for NNCP}

While there are multiple optimization methods to compute NNCP, we will focus on a class of methods that use Block Coordinate Descent (BCD), which is also known as the nonlinear Gauss-Seidel method \cite{Bertsekas1999}.
In BCD, the variables are partitioned into blocks, and each variable block is cyclically updated to optimality with all other blocks fixed.
For details on the convergence properties and comparisons of BCD methods for nonnegative matrix and tensor decomposition problems, see \cite{KHP2014}.
We consider BCD methods for NNCP that choose the entire factor matrices as the blocks, which is also often referred to as Alternating Least Squares.
In this case, every subproblem is a linear nonnegative least squares problem. Formally, the following problem 
is solved iteratively for $n = 1 \cdots N$:

$$\HH^{(n)} \leftarrow \argmin_{\HH \geq 0}  \left\| \Mn{S}{n} \HH^T - {\Mn{M}{n}}^T \right\|_F^2.$$
The BCD iteration is guaranteed to converge to a stationary point, but there is no guarantee that it converges to a global minimum.
The number of outer iterations to convergence is problem dependent but typically ranges from 10s to 1000s.


\Cref{alg:nncp} shows the pseudocode for BCD applied to NNCP.
\Cref{line:MTTKRP,line:GH,line:Gn} compute matrices involved in the gradients of the subproblem objective functions, and \cref{line:NLS} uses those matrices to update the current factor matrix.

The NLS-Update in \cref{line:NLS} can be implemented in different ways.
In a faithful BCD algorithm, the subproblems are solved exactly; in this case, the subproblem is a nonnegative linear least squares problem, which is convex.
We use the Block Principal Pivoting (BPP) method \cite{KP2011,KHP2014}, which is an active-set-like method, to solve the subproblem exactly.

However, as discussed in \cite{KBP2018} for the matrix case, there are other reasonable alternatives to updating the factor matrix without solving the subproblem exactly.
For example, we can more efficiently update individual columns of the factor matrix as is done in the Hierarchical Alternating Least Squares (HALS) method \cite{CP2009}.
In this case, the update rule is 
$$\Mn{H}{n}(:,r) \leftarrow \lt[ \Mn{H}{n}(:,r) + \Mn{M}{n}(:,r) - (\Mn{H}{n} \Mn{S}{n})(:,r)  \rt]_+$$
which involves the same matrices $\Mn{M}{n}$ and $\Mn{S}{n}$ as BPP.
Other possible BCD methods include Alternating Optimization and Alternating Direction Method of Multipliers (AO-ADMM) \cite{HSL2015, SBK2017} and Nestrov-based algorithms \cite{LKLHS2017}.
The parallel algorithm presented in this paper is generally agnostic to the approach used to solve the nonnegative least squares subproblems, as all these methods are bottlenecked by the subroutine they have in common, the MTTKRP.

\begin{algorithm}
\caption{$\CP = \text{NNCP}(\TA,R)$}
\label{alg:nncp}
\begin{algorithmic}[1]
\Require $\TA$ is $I_1\times \cdots \times I_N$ tensor, $R$ is approximation rank
\State \Comment{Initialize data}
\For{$n=2$ to $N$}
	\State Initialize $\Mn{H}{n}$ 
	\State $\Mn{G}{n} = \MnTra{H}{n}\Mn{H}{n}$
\EndFor
\State \Comment{Compute NNCP approximation}
\While{not converged}
	\State \Comment{Perform outer iteration of BCD}
	\For{$n=1$ to $N$}
	\State \Comment{Compute new factor matrix in $n$th mode}
	\State $\Mn{M}{n} = \text{MTTKRP}(\TA,\{\Mn{H}{i}\},n)$
		\label{line:MTTKRP}
	\State $\Mn{S}{n} = \Mn{G}{1} \Hada \cdots \Hada \Mn{G}{n-1} \Hada \Mn{G}{n+1} \Hada \cdots \Hada \Mn{G}{N}$
		\label{line:GH}
	\State $\Mn{H}{n} = \text{NLS-Update}(\Mn{S}{n},\Mn{M}{n})$
		\label{line:NLS}
	\State $\Mn{G}{n} = \MnTra{H}{n}\Mn{H}{n}$
		\label{line:Gn}
	\EndFor
\EndWhile
\Ensure $\TA \approx \CP$
\end{algorithmic}
\end{algorithm}

\subsection{Parallel Computation Model}

To analyze our algorithms we use the MPI model of distributed-memory parallel computation, where we assume a fully connected network. 
Sending a message of $W$ words from one processor to another costs $\alpha + \beta W$, where $\alpha$ is the latency and $\beta$ to be the per word or bandwidth cost. 
In particular, we will use collective communication over groups of $P$ processors, and we will assume the use of efficient algorithms \cite{TRG05,CH+07}.
In this case, an All-Reduce, which sums data initially distributed across processors and stores the result of size $W$ redundantly on every processor, costs $2\alpha \log P + 2\beta W (P-1)/P$.
An All-Gather collects data initially distributed across processors and stores the union of size $W$ redundantly on all processors and costs $\alpha \log P + \beta W (P-1)/P$.
A Reduce-Scatter sums data initially distributed across processors and partitions the result across processors, which costs $\alpha \log P + \beta W (P-1)/P$, where $W$ is the size of the data that each processor initially stores.
Reduction operations also include a flop cost but we will omit it because it is usually dominated by communication.

%% file: survey.tex

\section{Related Work} 
\label{sec:survey}

The formulation of NNCP with least squares error and algorithms for computing it go back to \cite{Paatero97,WW01}, developed in part as a generalization of nonnegative matrix factorization algorithms \cite{LS99} to tensors.
Sidiropoulos et al. \cite{SLFHPF2017} provide a more detailed and complete survey that includes basic tensor factorization models with and without constraints, broad coverage of algorithms, and recent driving applications.
The tensor operations discussed and the notation used in this paper follow Kolda and Bader's survey \cite{KB2009}. 

Recently, there has been growing interest in scaling tensor operations to bigger data and more processors in both the data mining/machine learning and the high performance computing communities. 
For sparse tensors, there have been parallelization efforts to compute CP decompositions both on shared-memory platforms \cite{SRSK2015,LCPSV17} as well as distributed-memory platforms \cite{KU16,SK16,KU18}, and these approaches can be generalized to constrained problems \cite{SBK2017}.
The focus of this work is on dense tensors, but many of the ideas for sparse tensors are applicable to the dense case, including parallel data distributions, communication pattern, and techniques to avoid recomputation across modes.

In particular, Liavas et al. \cite{LK+17b} extend a parallel algorithm designed for sparse tensors \cite{SK16} to the 3D dense case.
They use the ``medium-grained'' dense tensor distribution and row-wise factor matrix distribution, which is exactly the same as our distribution strategy (see \cref{sec:datadist}), and they use a Nesterov-based algorithm to enforce the nonnegativity constraints.
Their code is publicly available, and we compare our performance with theirs in \cref{sec:experiments}.
A similar data distribution and parallel algorithm for computing a single dense MTTKRP computation is proposed by Ballard, Knight, and Rouse \cite{BKR17-TR}. 
They prove that the algorithm is communication optimal, but they do not provide an implementation.
Another approach to parallelizing NNCP decomposition of dense tensors is presented by Phan and Cichocki \cite{PC11}, but they use a dynamic tensor factorization, which performs different, more independent computations across processors.

The idea of using dimension trees (discussed in \cref{sec:dimtrees}) to avoid recomputation within MTTKRPs across modes is introduced in \cite{PTC13a} for computing the CP decomposition of dense tensors.
It has also been used for sparse CP \cite{LCPSV17,KU18} and other tensor computations \cite{KU16}.

%
%


%% file: algorithm.tex

\section{Algorithm} 
\label{sec:algorithm}

\subsection{Dimension Trees}
\label{sec:dimtrees}

An important optimization of the CP-ALS algorithm is to re-use temporary values across inner iterations \cite{PTC13a,KU16-TR,LCPSV17,Kaya17}.
To illustrate the idea, consider a 3-way tensor $\T{X}$ approximated by $\dsquare{\M{U},\M{V},\M{W}}$ and the two MTTKRP computations $\Mn{M}{1}=\underline{\Mz{X}{1}(\M{W}}\Khat\M{V})$ and $\Mn{M}{2}=\underline{\Mz{X}{2}(\M{W}}\Khat\M{U})$ used to update factor matrices $\M{U}$ and $\M{V}$, respectively.
The underlined parts of the expressions correspond to the shared dependence of the outputs on the tensor $\T{X}$ and the third factor matrix $\M{W}$.
Indeed, a temporary quantity, which we refer to as a \emph{partial MTTKRP}, can be computed and re-used across the two MTTKRP expressions.
We refer to the computation that combines the temporary quantity with the other factor matrix to complete the MTTKRP computation as a multi-tensor-times-vector or \emph{multi-TTV}, as it consists of multiple operations that multiply a tensor times a set of vectors, each corresponding to a different mode.

To understand the steps of the partial MTTKRP and multi-TTV operations in more detail, we consider $\T{X}$ to be $I\times J\times K$ and $\M{U}$, $\M{V}$, and $\M{W}$ to have $R$ columns.
Then 
\begin{equation*}
\MnE{M}{1}{ir} = \sum_{i,j} \TE{X}{ijk} \ME{V}{jr} \ME{W}{kr} 
= \sum_{j} \ME{V}{jr} \sum_k \TE{X}{ijk} \ME{W}{kr} 
= \sum_{j} \ME{V}{jr} \TE{T}{ijr},
\end{equation*}
where $\T{T}$ is an $I\times J\times R$ tensor that is the result of a partial MTTKRP between tensor $\T{X}$ and the single factor matrix $W$.
Likewise,
\begin{equation*}
\MnE{M}{2}{jr} = \sum_{i,k} \TE{X}{ijk} \ME{U}{ir} \ME{W}{kr} 
= \sum_{i} \ME{U}{ir} \sum_k \TE{X}{ijk} \ME{W}{kr} 
= \sum_{i} \ME{U}{ir} \TE{T}{ijr},
\end{equation*}
and we see that the temporary tensor $\T{T}$ can be re-used.
From these expressions, we can also see that computing $\T{T}$ (a partial MTTKRP) corresponds to a matrix-matrix multiplication, and computing each of $\Mn{M}{1}$ and $\Mn{M}{2}$ from $\T{T}$ (a multi-TTV) corresponds to $R$ independent matrix-vector multiplications.
In this case, we compute $\Mn{M}{3}$ using a full MTTKRP.

For a larger number of modes, a more general approach can organize the temporary quantities to be used over a maximal number of MTTKRPs.
The general approach can yield significant benefit, decreasing the computation by a factor of approximately $N/2$ for dense $N$-way tensors.
The idea is introduced in \cite{PTC13a}, but we adopt the terminology and notation of \emph{dimension trees} used for sparse tensors in \cite{KU16-TR,Kaya17}.
In this notation, the root node is labeled $\{1,\dots,N\}$ and corresponds to the original tensor, a leaf is labeled $\{n\}$ and corresponds to the $n$th MTTKRP result, and an internal node is labeled by a set of modes $\{i,\dots,j\}$ and corresponds to a temporary tensor whose values contribute to the MTTKRP results of modes $i,\dots,j$.

\begin{figure}
\input{dimtree.tex}
\caption{Dimension tree example for $N=5$. 
The data associated with the root node is the original tensor, the data associated with the leaf nodes are MTTKRP results, and the data associated with internal nodes are temporary tensors.  
Edges labeled with PM correspond to partial MTTKRP computations, and edges labeled with mTTV correspond to multi-TTV computations.}
\label{fig:DT}
\end{figure}
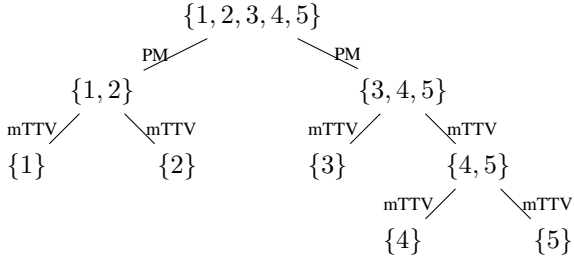

\Cref{fig:DT} illustrates a dimension tree for the case $N=5$.
Various shapes of binary trees are possible \cite{PTC13a,Kaya17}.
For dense tensors, the computational cost is dominated by the root's branches, which correspond to partial MTTKRP computations.
We perform the splitting of modes at the root so that modes are chosen contiguously with the respect to the layout of the tensor entries in memory.
In this way, each partial MTTKRP can be performed via BLAS's GEMM interface without reordering tensor entries in memory.
All other edges in a tree correspond to multi-TTVs and are typically much cheaper.
By organizing the memory layout of temporary quantities, the multi-TTV operations can be performed via a sequence of calls using BLAS's GEMV interface.
By using the BLAS in our implementation, we are able to obtain high performance and on-node parallelism.

\begin{figure}
\subfloat[Partial MTTKRP to compute node $\{3,4,5\}$ from root node $\{1,2,3,4,5\}$, executed via one GEMM call. \label{fig:DT-PM}]{\includegraphics[width=\columnwidth]{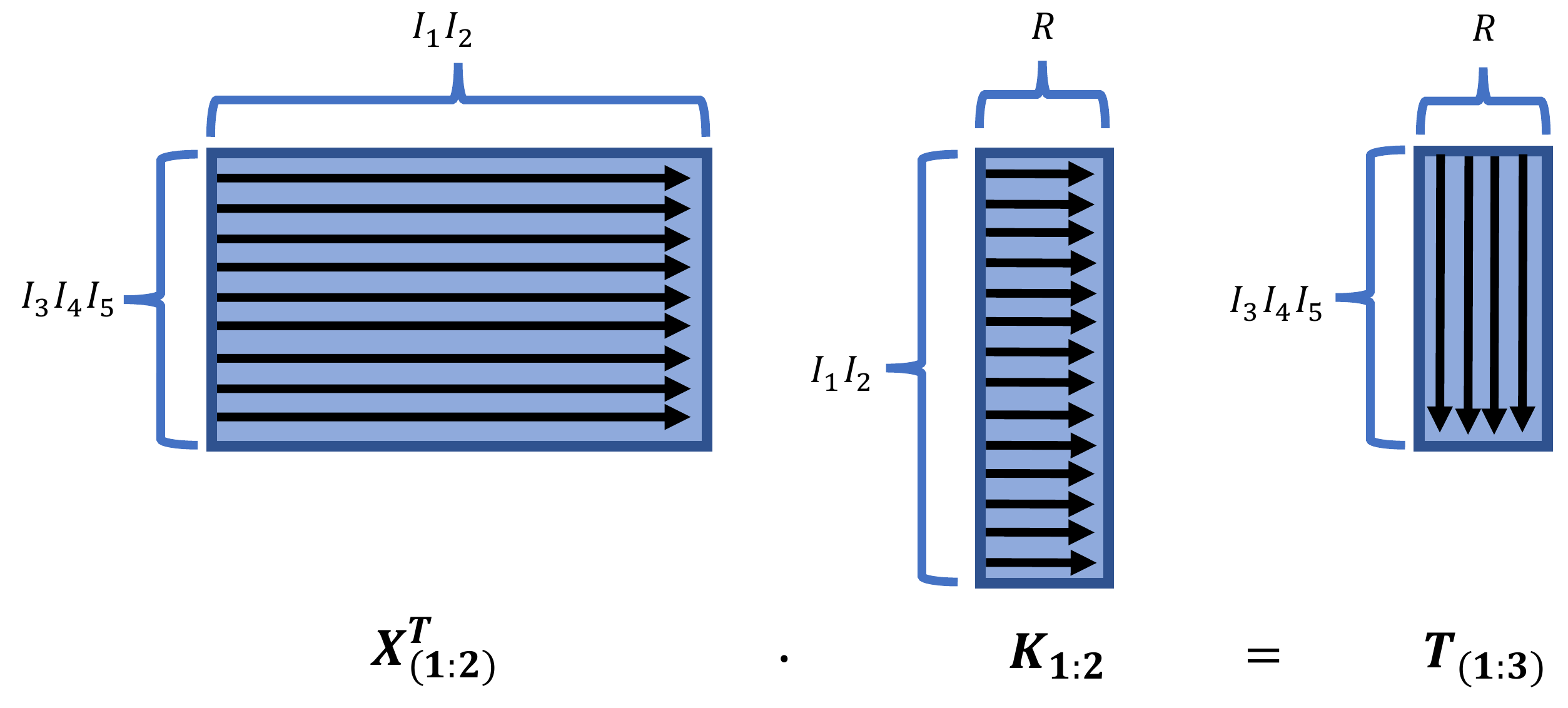}} \\
\subfloat[Multi-TTV to compute node $\{3\}$ from node $\{3,4,5\}$, executed via $R$ GEMV calls.  Here $\Mz{T}{1} \lbrack r\rbrack$ refers to the $r$th contiguous block of $\Mz{T}{1}$. \label{fig:DT-mTTV}]{\includegraphics[width=\columnwidth]{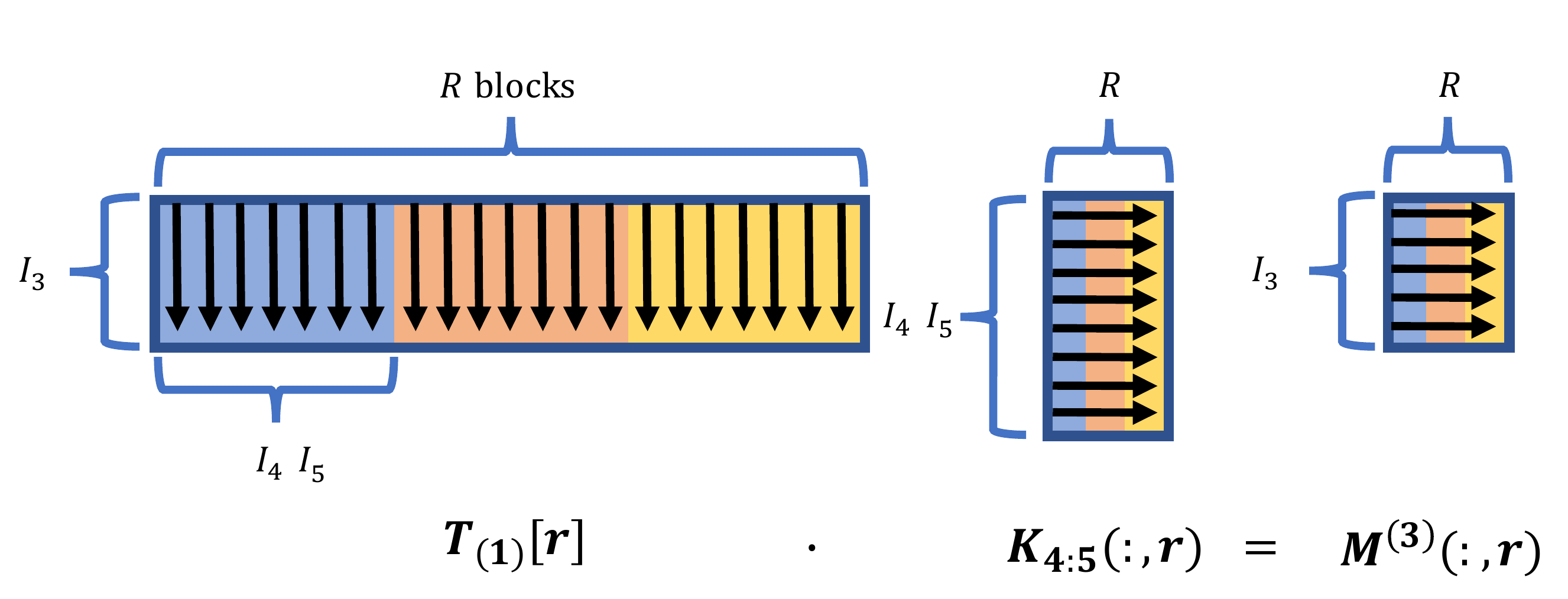}}
\caption{Data layout and dimensions for two example computations in dimension tree shown in \Cref{fig:DT}.  In this notation, $\Mz{X}{3:5}$ is the matricization of input tensor $\T{X}$ with respect to modes 3 through 5, $\M{K}_{1:2} = \Mn{H}{2} \Khat \Mn{H}{1}$, $\T T$ is the temporary $I_3 \times I_4 \times I_5 \times R$ tensor corresponding to node $\{3,4,5\}$ in the dimension tree, $\M{K}_{4:5} = \Mn{H}{5} \Khat \Mn{H}{4}$, and $\Mn{M}{3}$ is the MTTKRP result for mode 3.}
\label{fig:DTcartoon}
\end{figure}

\Cref{fig:DTcartoon} shows the data layout and dimensions of a partial MTTKRP and a multi-TTV taken from the example dimension tree in \Cref{fig:DT}.
\Cref{fig:DT-PM} shows a partial MTTKRP between the input tensor $\T{X}$ and the Khatri-Rao product of the factor matrices in modes 1 and 2, which produces a temporary tensor $\T{T}$ corresponding to the $\{3,4,5\}$ node in the dimension tree.
The key to efficiency in this computation is that the matricization of $\T{X}$ that assigns modes 1 through 2 to rows and modes 3 through 5 to columns is already column-major in memory.
Thus, we can use the GEMM interface and compute the temporary tensor $\T{T}$ without reordering any tensor entries.
\Cref{fig:DT-mTTV} depicts a multi-TTV that computes the results $\Mn{M}{3}$ from $\T{T}$ and the factor matrices in modes 4 and 5.
Here, the tensor $\T{T}$ is matricized with respect to only its first mode (of dimension $I_3$), but this matricization is also column-major in memory.
We choose the ordering of the modes of $\T{T}$ such that each of $R$ contiguous blocks is used to compute one column of the output matrix via a matrix-vector operation with a corresponding column of the Khati-Rao product of the other factor matrices.

No matter how the dimension tree is designed, the computational cost of each partial MTTKRP is $O(IR)$, where $I$ is the number of tensor entries and $R$ is the rank of the CP decomposition.
This is the same operation count as a full MTTKRP.
The computational cost of a multi-TTV is the number of entries in the temporary tensor, which is the product of a \emph{subset} of the original tensor dimensions multiplied by $R$.
Thus, it is computationally cheaper than the partial MTTKRPs, but it is also memory bandwidth bound.

The other subroutine necessary for implementing the dimension tree approach is the Khatri-Rao product of sets of factor matrices.
We implement the operation as a row-wise Hadamard product of a set of factor matrix rows, and we use OpenMP parallelization to obtain on-node parallelism.
The computational cost of this operation is also typically lower order, but the running time in practice suffers from also being memory bandwidth bound.

\subsection{Relative Error Computation}
\label{sec:error}

Given a model $\T{M}=\CP$, we compute the relative error $\|\TA - \T{M}\|/\|\TA\|$ efficiently by using the identity $\|\TA-\T{M}\|^2 = \|\TA\|^2 - 2\langle \TA, \T{M} \rangle + \|\T{M}\|^2.$
The quantity $\|\TA\|$ is fixed, and the other two terms can be computed cheaply given the temporary matrices computed during the course of the BCD algorithm.
The second term can be computed using the identity $\langle \TA, \T{M} \rangle = \langle \Mn{M}{N}, \Mn{H}{N} \rangle$, where $\Mn{M}{N} = \Mz{A}{N} (\Mn{H}{N-1} \Khat \cdots \Khat \Mn{H}{1})$ is the MTTKRP result in the $N$th mode.
The third term can be computed using the identity $\|\T{M}\|^2 = \V{1}^\Tra(\Mn{S}{N} \Hada \MnTra{H}{N} \Mn{H}{N})\V{1}$ where $\Mn{S}{N}=\MnTra{H}{1} \Mn{H}{1} \Hada \cdots \Hada \MnTra{H}{N-1} \Mn{H}{N-1}$.
Both matrices $\Mn{M}{N}$ and $\Mn{S}{N}$ are computed during the course of the BCD algorithm for updating the factor matrix $\Mn{H}{N}$.
The extra computation involved in computing the relative error is negligible.
These identities have been used previously \cite{KB09,TensorBox,SK16,LKLHS2017}.

\subsection{Parallel Algorithm}

\begin{algorithm}
\caption{$\CP = \text{Par-NNCP}(\TA,R)$}
\algfontsize
\label{alg:Par-NNCP-short}
\begin{algorithmic}[1]
\Require $\TA$ is an $I_1\times \cdots \times I_N$ tensor distributed across a $P_1\times \cdots \times P_N$ grid of $P$ processors, so that $\TA_{\V{p}}$ is $(I_1/P_1)\times \cdots \times (I_N/P_N)$ and is owned by processor $\V{p}=(p_1,\dots,p_N)$, $R$ is rank of approximation
\For{$n=2$ to $N$}
	\State Initialize $\Mn{H}{n}_{\V{p}}$ of dimensions $(I_n/P)\times R$ 
	\State $\M[\overline]{G} = \text{Local-SYRK}(\Mn{H}{n}_{\V{p}})$
	\State $\Mn{G}{n} = \text{All-Reduce}(\M[\overline]{G},\textsc{All-Procs})$
	\State $\Mn{H}{n}_{p_n} = \text{All-Gather}(\Mn{H}{n}_{\V{p}},\textsc{Proc-Slice}(n,\VE{p}{n}))$
\EndFor
\State \Comment{Compute NNCP approximation}
\While{not converged}
	\label{line:while}
	\State \Comment{Perform outer iteration of BCD}
	\For{$n=1$ to $N$}
		\label{line:for}
		\State \Comment{Compute new factor matrix in $n$th mode}
		\State $\M[\overline]{M} = \text{Local-MTTKRP}(\TA_{p_1\cdots p_N},\{\Mn{H}{i}_{p_i}\},n)$
			\label{line:locMTTKRP}
		\State $\Mn{M}{n}_{\V{p}} = \text{Reduce-Scatter}(\M[\overline]{M},\textsc{Proc-Slice}(n,\VE{p}{n}))$ 
			\label{line:reduce-scatter}
		\State $\Mn{S}{n} = \Mn{G}{1} \Hada \cdots \Hada \Mn{G}{n-1} \Hada \Mn{G}{n+1} \Hada \cdots \Hada \Mn{G}{N}$
			\label{line:hadamard}
		\State $\Mn{H}{n}_{\V{p}} = \text{Local-NLS-Update}(\Mn{S}{n},\Mn{M}{n}_{\V{p}})$
			\label{line:locNLS}
		\State \Comment{Organize data for later modes}
		\State $\M[\overline]{G} = {\Mn{H}{n}_{\V{p}}}^\Tra\Mn{H}{n}_{\V{p}}$
			\label{line:locSYRK}
		\State $\Mn{G}{n} = \text{All-Reduce}(\M[\overline]{G},\textsc{All-Procs})$
			\label{line:all-reduce}
		\State $\Mn{H}{n}_{p_n} = \text{All-Gather}(\Mn{H}{n}_{\V{p}},\textsc{Proc-Slice}(n,\VE{p}{n}))$
			\label{line:all-gather}
	\EndFor 
		\label{line:endfor}
\EndWhile
	\label{line:endwhile}
\Ensure $\TA \approx \CP$
\Ensure Local matrices: $\Mn{H}{n}_{\V{p}}$ is $(I_n/P)\times R$ and owned by processor $\V{p}=(p_1,\dots,p_N)$, for $1\leq n \leq N$, $\V{\lambda}$ stored redundantly on every processor
\end{algorithmic}
\end{algorithm}

\subsubsection{Algorithm Overview}

The basic sequential algorithm is given in \Cref{alg:nncp}, and the parallel version is given in \Cref{alg:Par-NNCP-short}.
We will refer to both the inner iteration, in which one factor matrix is updated (\cref{line:for} to \cref{line:endfor}), and the outer iteration, in which all factor matrices are updated (\cref{line:while} to \cref{line:endwhile}).
In the parallel algorithm, the processors are organized into a logical multidimensional grid (tensor) with as many modes as the data tensor.
The communication patterns used in the algorithm are all MPI collectives, including All-Reduce, Reduce-Scatter, and All-Gather.
The processor communicators (across which the collectives are performed) include the set of all processors and the sets of processors within the same processor slice.
Processors within a mode-$n$ slice all have the same $n$th coordinate.

The method of enforcing the nonnegativity constraints of the linear least squares solve (or update) generally affects only local computation because each row of a factor matrix can be updated independently.
In our algorithm, each processor solves the linear problem or computes the update for its subset of rows (see \cref{line:locNLS}). 
The most expensive (and most complicated) part of the parallel algorithm is the computation of the MTTKRP, which corresponds to \cref{line:locMTTKRP,line:reduce-scatter,line:all-gather}.

The details that are omitted from this presentation of the algorithm include the normalization of each factor matrix after it is computed and the computation of the residual error at the end of an outer iteration.
The computations do involve both local computation and communication, but their costs are negligible.
\ignore{A more detailed pseudocode is given in \Cref{alg:Par-NNCP-long}.}

\input{Algheader.tex}
\begin{figure*}[t]
\centering
  \subfloat[Start $n$th iteration with redundant subset of rows of each input matrix. \label{fig:inner_a}]{\input{AlgStart.tex}}  \quad
  \subfloat[Compute local MTTKRP for contribution to output matrix $\Mn{M}{2}$. \label{fig:inner_b}]{\input{AlgMTTKRP.tex}}  \quad
  \subfloat[Reduce-Scatter to compute and distribute rows of $\Mn{M}{2}$. \label{fig:inner_c}]{\input{AlgRS.tex}} \quad
  \subfloat[Compute local NLS update to obtain $\Mn{H}{2}_{\V{p}}$ from $\Mn{M}{2}_{\V{p}}$ (and $\Mn{S}{2}$). \label{fig:inner_d}]{\input{AlgNLS.tex}} \quad
  \subfloat[All-Gather to collect rows of $\Mn{H}{2}$ needed for later inner iterations. \label{fig:inner_e}]{\input{AlgAG.tex}}
  \caption{Illustration of 2nd inner iteration of Par-NNCP algorithm for 3-way tensor on a $3\times3\times3$ processor grid, showing data distribution, communication, and computation across steps.  Highlighted areas correspond to processor $(1,3,1)$ and its processor slice with which it communicates.  The column normalization and computation of $\Mn{G}{2}$, which involve communication across all processors, is not shown here.}
  \label{fig:inner} 
\end{figure*}
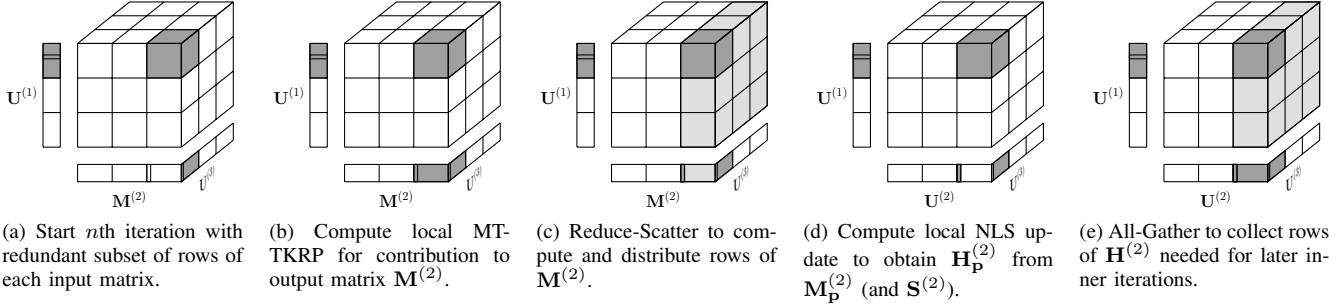

\subsubsection{Data Distribution}
\label{sec:datadist}

Given a logical processor grid of processors $P_1\times \cdots \times P_N$, we distribute the tensor $\TA$ in a block or Cartesian partition.
Each processor owns a local tensor of dimensions $(I_1/P_1)\times \cdots \times (I_N/P_N)$, and only one copy of the tensor is stored.
Locally, the tensor is stored linearly, with entries ordered in a natural mode-descending way that generalizes column-major layout of matrices.
Given a processor $\V{p}=(\VE{p}{1},\dots,\VE{p}{N})$, we denote its local tensor $\TA_{\V{p}}$.

Each factor matrix is distributed across processors in a block row partition, so that each processor owns a subset of the rows.
We use the notation $\Mn{H}{n}_{\V{p}}$, which has dimensions $I_n/P\times R$ to denote the local part of the $n$th factor matrix stored on processor $\V{p}$.
However, we also make use a redundant distribution of the factor matrices across processors, because all processors in a mode-$n$ processor slice need access to the same entries of $\Mn{H}{n}$ to perform their computations.
The notation $\Mn{H}{n}_{\VE{p}{n}}$ denotes the $I_n/P_n\times R$ submatrix of $\Mn{H}{n}$ that is redundantly stored on all processors whose $n$th coordinate is $\VE{p}{n}$ (there are $P/P_n$ such processors).

Other matrices involved in the algorithm include $\Mn{M}{n}_{\V{p}}$, which is the result of the MTTKRP computation and has the same distribution scheme as $\Mn{H}{n}_{\V{p}}$, and $\Mn{G}{n}$, which is the $R\times R$ Gram matrix of the factor matrix $\Mn{H}{n}$ and is stored redundantly on all processors.

\subsubsection{Inner Iteration}

The inner iteration is displayed graphically in \Cref{fig:inner} for a 3-way example and an update of the $2$nd factor matrix.
The main idea is that at the start of the $n$th inner iteration (\Cref{fig:inner_a}), all of the data is in place for each processor to perform a local MTTKRP computation.
This means that all processors in a slice redundantly own the same rows of the corresponding factor matrix (for all modes except $n$).
After the local MTTKRP is computed (\Cref{fig:inner_b}), each processor has computed a contribution to a subset of the rows of the global MTTKRP $\Mn{M}{n}$, but its contribution must be summed up with the contributions of all other processors in its mode-$n$ slice.
This summation is performed with a Reduce-Scatter collective across the mode-$n$ processor slice that achieves a row-wise partition of the result (in \Cref{fig:inner_c}, the light gray shading corresponds to the rows of $\Mn{M}{2}$ to which processor $(1,3,1)$ contributes and the dark gray shading corresponds to the rows it receives as output).
The output distribution of the Reduce-Scatter is designed so that afterwards, the update of the factor matrix in that mode can be performed row-wise in parallel.
Along with $\Mn{S}{n}$, which can be computed locally, each processor updates its own rows of the factor matrix given its rows of the MTTKRP result (\Cref{fig:inner_d}).
The remainder of the inner iteration is preparing and distributing the new factor matrix data for future inner iterations, which includes an All-Gather of the newly computed factor matrix $\Mn{H}{n}$ across mode-$n$ processor slices (\Cref{fig:inner_e}) and recomputing $\Mn{G}{n}={\Mn{H}{n}}^\Tra\Mn{H}{n}$.

\subsubsection{Analysis}

\begin{table*}
\centering
\begin{tabular}{|c|cccc|}
\hline
 & \textbf{Computation} & \textbf{Memory - Local MTTKRP} & \textbf{Communication} & \textbf{Memory - Par.~Algorithm}  \\
\hline
Par-NNCP w/ opt.~proc.~grid & $\frac{IR}{P}$ & $\frac{RI^{1/2}}{P^{1/2}}$ & $\frac{NRI^{1/N}}{P^{1/N}}$ & $\frac{NRI^{1/N}}{P^{1/N}}$  \\
Par-NNCP w/ gen.~proc.~grid & $\frac{IR}{P}$ & $\frac{RI^{1/2}}{P^{1/2}}$ & $\mathbf{R\sum_n \frac{I_n}{P_n}}$ & $\mathbf{R\sum_n \frac{I_n}{P_n}}$  \\
Par-NNCP w/o dim.~tree opt. & $\mathbf{\frac{NIR}{P}}$ & $\displaystyle \mathbf{\max_n \frac{RI/I_n}{P/P_n}}$  &  $\frac{NRI^{1/N}}{P^{1/N}}$ & $\frac{NRI^{1/N}}{P^{1/N}}$ \\
\hline
\end{tabular}
\smallskip
\caption{Leading-order per-iteration costs in terms of computation (flops), communication (words moved), and memory (words).  We ignore constants but omit big-Oh notation for clarity.  The first line corresponds to \cref{alg:Par-NNCP-short} with an optimal choice of processor grid and applying the dimension tree optimization locally.  The second line corresponds to the changes in communication and memory for a general processor grid.  The third line corresponds to the changes in computation and memory if the dimension tree optimization is not applied.}
\label{tab:costs}
\end{table*}

We will analyze the cost of a single outer iteration.
While the number of outer iterations is sensitive to the NLS method used, the outer iteration time is generally the same across NLS methods.
We summarize the analysis in \Cref{tab:costs}, showing the differences with and without an optimal processor grid and with and without using a dimension tree.

\paragraph{\emph{Computation}}
The local computation occurs at \cref{line:locMTTKRP,line:hadamard,line:locNLS,line:locSYRK}.
The cost of \cref{line:hadamard} is $O(NR^2)$, the cost of \cref{line:locNLS} is $O(R^3I_n/P)$, which is a loose upper bound for BPP and other methods \cite{KBP16}, and the cost of \cref{line:locSYRK} is $O(R(I_n/P)^2)$.
The sum of these three costs across all inner iterations is $O(R^2N^2+(R^3/P)\sum I_n+(R/P^2)\sum I_n^2)$, which is dominated by the cost of the MTTKRP.
When using dimension trees to perform the MTTKRP (\cref{line:locMTTKRP}), we compute the cost amortized over all inner iterations.
In this case, the cost is dominated by the two partial MTTKRP computations (from the root of the tree), which together are $O((R/P) \prod I_n)=O(IR/P)$ and dominate the costs of the multi-TTVs.
We note that this cost involves the product of all the tensor dimensions, which is why it dominates, and we note that it scales linearly with $P$.

\paragraph{\emph{Communication}}
The communication within the inner iteration occurs at \cref{line:reduce-scatter,line:all-reduce,line:all-gather}.
\Cref{line:all-reduce} involves $O(R^2)$ data and a collective across all processors.
\Cref{line:reduce-scatter,line:all-gather} involve $O(I_nR/P_n)$ data across a subset of $P/P_n$ processors.
Thus, the All-Reduce dominates the latency cost and the Reduce-Scatter/All-Gather dominate the bandwidth cost, for a total outer iteration communication cost of $O(R\sum I_n/P_n)$ words and $O(N\log P)$ messages.
If the optimal processor grid can be chosen to minimize communication (assuming $P$ is sufficiently factorable), then the bandwidth cost can achieve a value of $O(NRI^{1/N}/P^{1/N})$ by making the local tensors as cubical as possible.
We note that this cost scales with $P^{1/N}$, which is far from linear scaling.
However, it is proportional to the geometric mean of the tensor dimensions (on the order of one tensor dimension), which is much less than the computation cost dependence on the product of all dimensions.

\paragraph{\emph{Memory}}
The algorithm requires extra local memory to run.
Aside from the memory required to store the local tensor of $O(I/P)$ words and factor matrices of cumulative size $O((R/P)\sum I_n)$, each processor must be able to store a redundant subset of the rows of the factor matrices it needs to perform MTTKRP computations.
This corresponds to storing $P/P_n$ redundant copies of every factor matrix, which results in a local memory requirement of $O(R \sum I_n/P_n)$ for a general processor grid.
The processor grid that minimizes communication also minimizes local memory, and the extra memory requirement can be as low as $O(NRI^{1/N}/P^{1/N})$, which is typically dominated by $O(I/P)$.

The dimension tree algorithm also requires extra temporary memory space, but the space required tends to be much smaller than what is required to store the local tensor.
If the tensor dimensions can be partitioned into two parts with approximately equal geometric means, the extra memory requirement for running a dimension tree is as small as $O(R\sqrt{I/P})$, which is also typically dominated by $O(I/P)$.

%% file: dimtree.tex

\begin{center}
\begin{tikzpicture}

\node (12345) at (4,2) {$\{1,2,3,4,5\}$};
\node (12) at (2,1) {$\{1,2\}$};
\node (345) at (6,1) {$\{3,4,5\}$};
\node (1) at (1,0) {$\{1\}$};
\node (2) at (3,0) {$\{2\}$};
\node (3) at (5,0) {$\{3\}$};
\node (45) at (7,0) {$\{4,5\}$};
\node (4) at (6,-1) {$\{4\}$};
\node (5) at (8,-1) {$\{5\}$};

\scriptsize
\path[draw] (12345) edge [left,align=right] node {PM} (12);
\path[draw] (12345) edge [right,align=left] node {PM} (345);
\path[draw] (12) edge [left,align=right] node {mTTV} (1);
\path[draw] (12) edge [right,align=left] node {mTTV} (2);
\path[draw] (345) edge [left,align=right] node {mTTV} (3);
\path[draw] (345) edge [right,align=left] node {mTTV} (45);
\path[draw] (45) edge [left,align=right] node {mTTV} (4);
\path[draw] (45) edge [right,align=left] node {mTTV} (5);
\normalsize

\end{tikzpicture}
\end{center}

%% file: Algheader.tex

\newcommand{\procdim}{3}
\newcommand{\proc}{\draw[black,shift={(-.5,-.5)}] (0,0) grid (\procdim,\procdim);}
\newcommand{\highlight}{gray!75}
\newcommand{\commhighlight}{gray!25}
\newcommand{\parscale}{.46}
\newcommand{\secfaclabel}{$\Mn{M}{2}$}

\newcommand{\parbasepic}{
\begin{scope}[canvas is yz plane at x=.5,shift={(1.5,-1.5)}]
	\draw[fill=\highlight,shift={(0,1)}] (0,0) rectangle (1,1);
	\draw[fill=\highlight,shift={(-2.5,2-1/3)},xscale=.5] (0,0) rectangle (-1,-1/9);
	\draw[shift={(0,-1.5)},yscale=.5] (0,0) rectangle (1/9,-1);
\end{scope}
\begin{scope}[canvas is zx plane at y=(\procdim-.5),rotate=-90,shift={(-.5,-3)}]
	\draw[fill=\highlight,shift={(0,2.5)}] (0,0) rectangle (1,1);
	\draw[fill=\highlight,yscale=.5] (0,0) rectangle (1/9,-1);
\end{scope}
\begin{scope}[canvas is yx plane at z=.5,yscale=-1,rotate=0]
	\draw[fill=\highlight,shift={(1.5,-.5)}] (0,0) rectangle (1,1);
\end{scope}

\begin{scope}[canvas is yz plane at x=.5,rotate=-90]
	\proc
\end{scope}
\begin{scope}[canvas is yx plane at z=.5,yscale=-1,rotate=0]
	\proc
\end{scope}
\begin{scope}[canvas is zx plane at y=(\procdim-.5),rotate=180]
	\proc
\end{scope}

\begin{scope}[canvas is yz plane at x=.5,shift={(1.5,-.5)}]
	\draw[shift={(-2.5,0)},xscale=.5] (0,-2) grid (-1,1);
	\node[draw=none] at (-3.6,-.5) {\Large $\Mn{U}{1}$};
	\draw[shift={(0,-2.5)},yscale=.5] (-2,0) grid (1,-1);
	\node[draw=none] at (-.5,-3.5) {\Large \secfaclabel};
\end{scope}
\begin{scope}[canvas is zx plane at y=(\procdim-.5),rotate=-90,shift={(1.5,-.5)}]
	\draw[shift={(0,-2.5)},yscale=.5] (-2,0) grid (1,-1);
	\node[draw=none] at (-.5,-3.5) {\Large $\Mn{U}{3}$};
\end{scope}
}

%% file: AlgStart.tex

\begin{tikzpicture}[x={(-0.5cm,-0.4cm)}, y={(1cm,0cm)}, z={(0cm,1cm)},every node/.append style={transform shape},scale=\parscale]


\begin{scope}[canvas is yz plane at x=.5]
	\draw[fill=\highlight,shift={(-1,.5)},xscale=.5] (0,0) rectangle (-1,-1);
\end{scope}
\begin{scope}[canvas is zx plane at y=(\procdim-.5),rotate=-90]
	\draw[fill=\highlight,yscale=.5,shift={(-.5,-6)}] (0,0) rectangle (1,-1);
\end{scope}

\begin{scope}[canvas is yz plane at x=.5]
\end{scope}

\parbasepic

\end{tikzpicture}

%% file: AlgMTTKRP.tex

\begin{tikzpicture}[x={(-0.5cm,-0.4cm)}, y={(1cm,0cm)}, z={(0cm,1cm)},every node/.append style={transform shape},scale=\parscale]


\begin{scope}[canvas is yz plane at x=.5]
	\draw[fill=\highlight,shift={(-1,.5)},xscale=.5] (0,0) rectangle (-1,-1);
\end{scope}
\begin{scope}[canvas is zx plane at y=(\procdim-.5),rotate=-90]
	\draw[fill=\highlight,yscale=.5,shift={(-.5,-6)}] (0,0) rectangle (1,-1);
\end{scope}

\begin{scope}[canvas is yz plane at x=.5]
	\draw[fill=\highlight,shift={(1.5,-3)},yscale=.5] (0,0) rectangle (1,-1);
\end{scope}

\parbasepic

\end{tikzpicture}

%% file: AlgRS.tex

\begin{tikzpicture}[x={(-0.5cm,-0.4cm)}, y={(1cm,0cm)}, z={(0cm,1cm)},every node/.append style={transform shape},scale=\parscale]


\begin{scope}[canvas is yz plane at x=.5]
	\draw[fill=\commhighlight,shift={(1.5,.5)}] (0,0) rectangle (1,-3);
	\draw[fill=\commhighlight,shift={(1.5,-3)},yscale=.5] (0,0) rectangle (1,-1);
	\draw[fill=\highlight,shift={(1.5,-3)},yscale=.5] (0,0) rectangle (1/9,-1);
	\draw[fill=\highlight,shift={(-1,.5)},xscale=.5] (0,0) rectangle (-1,-1);
\end{scope}
\begin{scope}[canvas is zx plane at y=(\procdim-.5),rotate=-90]
	\draw[fill=\commhighlight,shift={(-.5,.5)}] (0,0) rectangle (3,-3);
	\draw[fill=\highlight,yscale=.5,shift={(-.5,-6)}] (0,0) rectangle (1,-1);
\end{scope}
\begin{scope}[canvas is yx plane at z=.5,yscale=-1,rotate=0]
	\draw[fill=\commhighlight,shift={(1.5,-.5)}] (0,0) rectangle (1,3);
\end{scope}

\parbasepic

\end{tikzpicture}

%% file: AlgNLS.tex

\renewcommand{\secfaclabel}{$\Mn{U}{2}$}

\begin{tikzpicture}[x={(-0.5cm,-0.4cm)}, y={(1cm,0cm)}, z={(0cm,1cm)},every node/.append style={transform shape},scale=\parscale]


\begin{scope}[canvas is yz plane at x=.5]
	\draw[fill=\highlight,shift={(1.5,-3)},yscale=.5] (0,0) rectangle (1/9,-1);
	\draw[fill=\highlight,shift={(-1,.5)},xscale=.5] (0,0) rectangle (-1,-1);
\end{scope}
\begin{scope}[canvas is zx plane at y=(\procdim-.5),rotate=-90]
	\draw[fill=\highlight,yscale=.5,shift={(-.5,-6)}] (0,0) rectangle (1,-1);
\end{scope}
\begin{scope}[canvas is yx plane at z=.5,yscale=-1,rotate=0]
\end{scope}

\parbasepic

\end{tikzpicture}

%% file: AlgAG.tex

\renewcommand{\secfaclabel}{$\Mn{U}{2}$}

\begin{tikzpicture}[x={(-0.5cm,-0.4cm)}, y={(1cm,0cm)}, z={(0cm,1cm)},every node/.append style={transform shape},scale=\parscale]


\begin{scope}[canvas is yz plane at x=.5]
	\draw[fill=\commhighlight,shift={(1.5,.5)}] (0,0) rectangle (1,-3);
	\draw[fill=\highlight,shift={(1.5,-3)},yscale=.5] (0,0) rectangle (1,-1);
	\draw[fill=\highlight,shift={(1.5,-3)},yscale=.5] (0,0) rectangle (1/9,-1);
	\draw[fill=\highlight,shift={(-1,.5)},xscale=.5] (0,0) rectangle (-1,-1);
\end{scope}
\begin{scope}[canvas is zx plane at y=(\procdim-.5),rotate=-90]
	\draw[fill=\commhighlight,shift={(-.5,.5)}] (0,0) rectangle (3,-3);
	\draw[fill=\highlight,yscale=.5,shift={(-.5,-6)}] (0,0) rectangle (1,-1);
\end{scope}
\begin{scope}[canvas is yx plane at z=.5,yscale=-1,rotate=0]
	\draw[fill=\commhighlight,shift={(1.5,-.5)}] (0,0) rectangle (1,3);
\end{scope}

\parbasepic

\end{tikzpicture}

%% file: experiments.tex

\section{Performance Results} \label{sec:experiments}

\input{plots}

\subsection{Datasets}

\subsubsection{Hyperspectral Images (HSI)}
For comparison with previous work \cite{LK+17b}, we consider the same 3D hyperspectral imaging dataset called ``Souto\_wood\_pile'' \cite{FAN16}. 
NNCP is often used on HSI data sets for classification and blind source separation of materials with differing spectral signatures.
The hyperspectral datacube has dimensions $1024 \times 1344 \times 33$ and represents a set of 33 grayscale images of size 1344 $\times$ 1024 pixels sampled at wavelengths 400, 410, $\dots$, 720 nm, with each pixel value representing spectral radiance in $W m^{-2} sr^{-1} nm^{-1}$. 
We also consider the Nogueir\'{o} scene dataset, which is a sequence of 9 time-lapse HSI images of the same scene acquired at about 1-hour intervals.
In each scene, hyperspectral images were acquired at about 1-hour intervals. 
Each Nogueir\'{o} scene HSI image has the same properties as the Souto\_wood\_pile data set, so the corresponding tensor has dimensions $1024 \times 1344 \times 33 \times 9$. 

\subsubsection{Dynamic Functional Connectivity (dFC)}
We also consider dynamic functional connectivity datasets that are generated from fMRI images of human brains.
Given a 4D fMRI data set of voxel measurements across multiple timesteps, voxels containing brain data are partitioned into a set of regions of interest (specified using domain-specific knowledge), and a single time-series signal is aggregated for each region of interest.
Then, an instantaneous correlation is computed for each time point and pair of regions, and this process is repeated for a number of subjects.
Computing a CP decomposition of this data helps to discover patterns of brain connectivity among different regions and also differentiate among individuals.
For our representative dFC data set, we consider 246 brain regions, which yields 30{,}012 unique pairs of regions, 1200 times steps, and 500 subjects, or a tensor of dimension $30{,}012\times 1200 \times 500$ \cite{VEU+12,THBGW17}.

\subsubsection{Synthetic}
Our synthetic data sets are constructed from a CP model with an exact low rank with no added noise.
In this case we can confirm that the residual error of our algorithm with a random start converges to zero.
For the purposes of benchmarking, we run a fixed number of iterations of the BCD algorithm rather than using a convergence check.


\subsection{Machine Details}
The entire experimentation was performed on Eos, a supercomputer at the Oak Ridge Leadership Computing Facility. 
Eos is a 736-node Cray XC30 cluster of Intel Xeon E5-2670 processors with a total of 47.104TB of memory. 
Its compute nodes are organized in blades where each blade contains 4 nodes, and every node has 2 sockets with 8 physical cores and 64GB memory. 
The machine support Intel's hyper-threading (HT), but we restricted it because HT offers minimal improvement for BLAS and LAPACK operations. 
In total, the Eos compute partition contains 11,776 traditional processor cores and our experiments used up to 4,096 cores (35\% of the machine). 

Our objective of the implementation is using open source software as much as possible 
to promote reproducibility and reuse of our code.
We use Armadillo \cite{sanderson2010} for matrix representation
and operations.  
In Armadillo, the elements of the dense matrix are stored in column major order.
For dense BLAS and LAPACK operations, we linked Armadillo with the default LAPACK/BLAS wrappers from Cray. 
For compiler, we use GNU C++ Compiler (g++ (GCC) 6.3.0) and MPI library is from Cray.  We could also 
compile and run the code in Rhea the commodity cluster from OLCF with entire open source libraries such as OpenBLAS and OpenMPI. 

\subsection{Comparison Implementations}
The implementation proposed by Liavas et al. \cite{LK+17b} is the only publicly available distributed-memory software (of which we are aware) for computing the CP decomposition of dense tensors, with or without constraints.
We use the acronym NbAO-NTF for Nesterov-based Alternating Optimization Nonnegative Tensor Factorization to refer to their code.

It is based on the same parallel algorithm as our implementation, though it is limited to 3D tensors.
The code uses MPI collectives for communication and Eigen \cite{eigenweb} as an interface to BLAS and LAPACK.
We compiled the code linked to BLAS/LAPACK wrappers from Cray  (the same BLAS implementation used by our code) but we were unable to run multithreaded BLAS with their code.
For fair comparison, we use a flat MPI configuration (one MPI process per core) on all comparisons between the two implementations.

We also point out a difference between the Nesterov-based algorithm and the BPP algorithm for solving the NLS subproblems.
The Nesterov-based algorithm attempts an acceleration step using a linear combination of the current and proposed future step; however, it re-computes the residual error before deciding whether or not to accept or reject the acceleration step.
This residual error cannot always be computed cheaply, using the technique described in \cref{sec:error}, and it contributes significantly (approximately 25\%) to the overall run time.
Because the BPP algorithm does not require this extra computation, and studying convergence behavior of the different NLS algorithms is beyond the scope of this work, we remove the time spent in the acceleration step of NbAO-NTF in all our comparisons.

Our proposed algorithm uses dimension trees, but we also benchmark our implementation without that optimization to highlight its importance.
We use an existing implementation to perform the individual MTTKRPs \cite{HBJT18} with this approach.

\subsection{Strong Scaling}

We perform two strong scaling experiments to compare performance with NbAO-NTF.
The experiments use a cubical synthetic tensor and the HSI image used in \cite{LK+17b}, both of which are 3D.

\begin{figure}
\begin{tikzpicture}
\renewcommand{\datafile}{str_3D_syn.dat}
\renewcommand{\numiterations}{1}
\liavastrue
\strongscalingplot
\end{tikzpicture}
\caption{Strong scaling of 3D synthetic tensor with dimension $1024\times 1024\times 1024$ on processor grids $2^k\times 2^k\times 2^k$ for $k\in\{0,1,2,3,4\}$.  The rank is fixed at 32.}
\label{fig:strongsynthetic3D}
\end{figure}
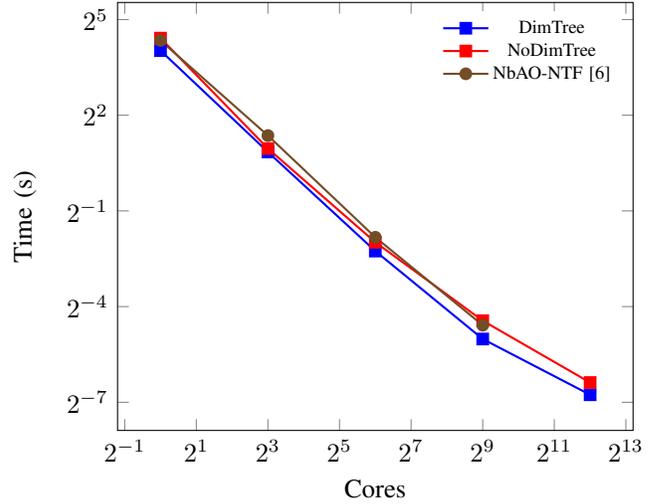

The performance on the cubical synthetic tensor is shown in \Cref{fig:strongsynthetic3D}. 
We can observe from the figure that all the three algorithms scale nearly linearly as the problem remains compute bound. 
Our algorithm (with the dimension tree optimization) achieves a speedup of $1771\times$ on 4096 cores over the same implementation running on 1 core.
Recall from \label{sec:analysis} that the computation scales linearly with $1/P$ while the communication scales with $1/P^{1/N}=1/P^{1/3}$. 
As is evident from the figure, the communication cost does not degrade performance even for thousands of cores. 
Our proposed algorithm with dimension trees is 35\% faster than NbAO-NTF at 512 cores (with similar relative difference for other core counts).
This performance improvement is due in large part to the 50\% reduction in arithmetic operations provided by the dimension tree optimization.
There is little difference in performance between our implementation without dimension trees and NbAO-NTF.

\begin{figure}
\begin{tikzpicture}
\renewcommand{\datafile}{str_3D_HSI.dat}
\renewcommand{\numiterations}{10}
\liavastrue
\strongscalingplot
\end{tikzpicture}
\caption{Strong scaling of 3D HSI real world data with dimension 1024 x1344 x 33 on processor grids of $k \times k\times 1$ for $k \in {1, 2, 4, 8, 16, 32}$. The rank is fixed at 32.}
\label{fig:stronghsi3D}
\end{figure}
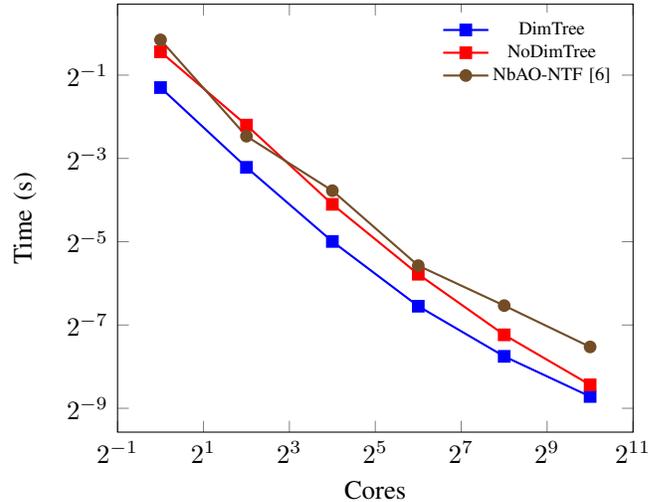

\Cref{fig:stronghsi3D} shows the strong scaling on the HSI data. 
In this case, our proposed algorithm with dimension trees is over $2\times$ faster than NbAO-NTF, but part of this speedup is due to differences in the NLS update algorithms.
For the low core count, the dimension tree provides a 60\% speedup compared to the MTTKRP time in NbAO-NTF.
At the high core counts for this experiment, the local MTTKRP is no longer the dominating cost.


\begin{figure}
\begin{tikzpicture}
\renewcommand{\datafile}{str_5D_syn.dat}
\renewcommand{\numiterations}{10}
\liavasfalse
\strongscalingplot
\end{tikzpicture}
\caption{Strong scaling of 5D synthetic tensor with dimension $64\times 64\times 64\times 64\times 64$ on processor grids $1\times1\times1\times1\times1$, $2\times1\times1\times1\times1$, $\dots$, $2\times2\times2\times2\times2$.  The rank is fixed at 32.}
\label{fig:strongsynthetic5D}
\end{figure}
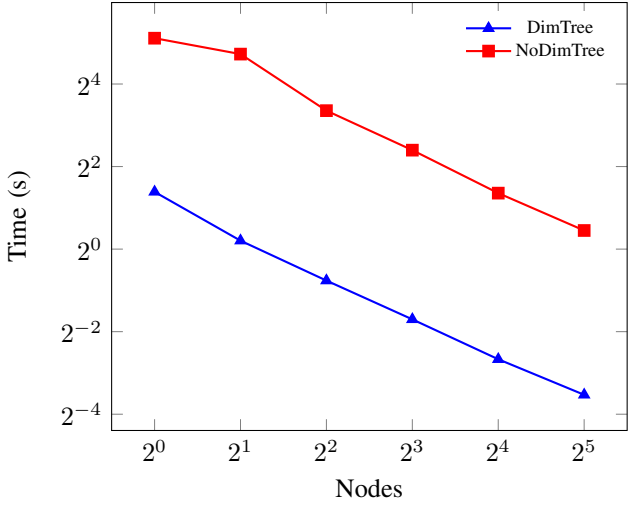

In \Cref{fig:strongsynthetic5D}, we benchmark performance for a 5D cubical tensor with each dimension set to 64.
Because the tensor is 5D, we can no longer compare against NbAO-NTF.
We see a $13-16\times$ speed up using a dimension tree over not using one for this problem.
As predicted, the dimension tree optimization saves relatively more arithmetic for higher-order tensors.
However, the reduction in leading order flop cost is only $2.5\times$ for $N=5$; the rest of the speedup comes from more efficient DGEMM performance and avoiding memory-bound KRP computation.  
That is, although the flop count of KRP computation is lower order, it still contributes to the run time because it is inefficient.
Also, for tensors with balanced dimensions, the dimension tree approach yields more favorable shapes for DGEMM.

\subsection{Weak Scaling Time Breakdown}

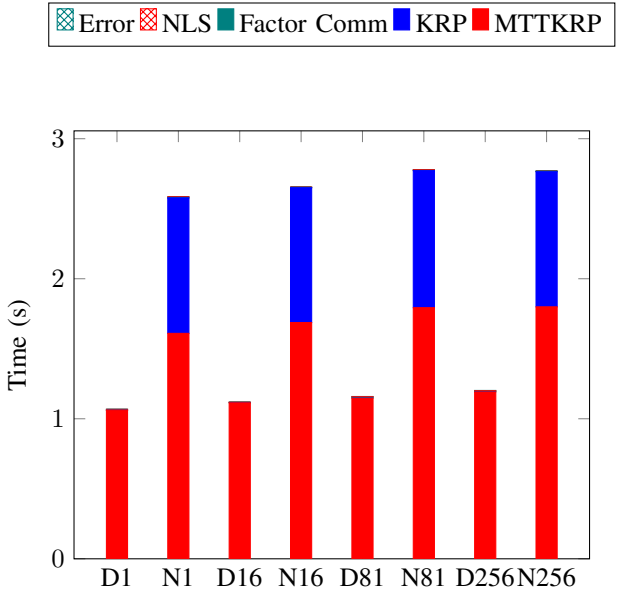
\begin{figure}
\begin{tikzpicture}
\renewcommand{\datafile}{wk_4D_syn.dat}
\renewcommand{\numiterations}{10}
\weakthreeDfalse
\breakdownplot
\end{tikzpicture}
\caption{Weak scaling of 4D synthetic tensors with (D) and without (N) the use of dimension trees.  The tensor and processor grid dimensions are $128k\times 128k\times 128k\times 128k$ and $k\times k\times k\times k$ for $k\in\{1,2,3,4\}$, and the rank is fixed at 32.  The reported times are per iteration.}
\label{fig:weaksynthetic4D}
\end{figure}

We also perform a weak scaling experiment to understand the time it takes to solve bigger problems with more processors.
In this experiment, we use a synthetic 4D tensor and keep the amount of tensor data assigned to each processor constant, with tensor and processor grid of dimensions $128k\times 128k\times 128k\times 128k$ and $k\times k\times k\times k$ for $k\in\{1,2,3,4\}$ and the rank fixed at 32. 
The results of the breakdown plot is shown in \Cref{fig:weaksynthetic4D}. 
In this case, the algorithm is compute bound with and without the use of the dimension tree, so the total time of the weak scaling remains fixed for both cases. 
Using a dimension tree, the time is completely dominated by the MTTKRP computation.  
Without using a dimension tree, we observe that the KRP is expensive and yields a $2.5\times$ slower total run time even in the 4D case. 

\subsection{Varying Processor Grid}

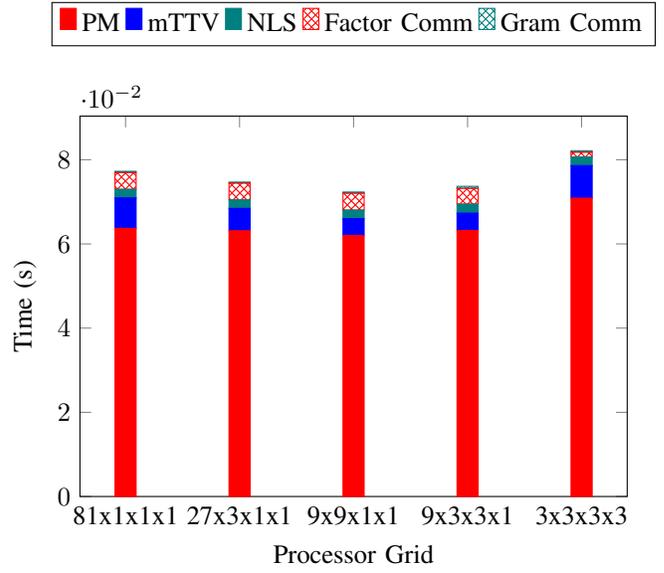
\begin{figure}
\renewcommand{\datafile}{comm_sweep_scale_pgf_10.dat}
\renewcommand{\numiterations}{10}
\begin{tikzpicture}
\begin{axis}[	
	ybar stacked,
	bar width=8pt,
	width=\columnwidth,
	height =.75\columnwidth,
	ylabel={Time (s)}, 
	xlabel={Processor Grid},
	y label style={yshift=-.5cm},
	ymin=0,
	symbolic x coords={81x1x1x1-10, 27x3x1x1-10,9x9x1x1-10,9x3x3x1-10,3x3x3x3-10,},
	xticklabels={,81x1x1x1,27x3x1x1,9x9x1x1,9x3x3x1,3x3x3x3},
	legend style={at={(0.5,1.3)},anchor=north},
	legend columns=-1,
]
	\setcolors
	\addplot table[x=algo-k, y expr=(\thisrow{mttkrp}/(\minvalue*\numiterations))] {\datafile};
	\addplot table[x=algo-k, y expr=(\thisrow{mttv}/(\minvalue*\numiterations))] {\datafile};
	\addplot table[x=algo-k, y expr=((\thisrow{nnls}+\thisrow{gram})/(\minvalue*\numiterations))] {\datafile};
	\addplot table[x=algo-k, y expr=((\thisrow{allgather}+\thisrow{reducescatter})/(\minvalue*\numiterations))] {\datafile};
	\addplot table[x=algo-k, y expr=(\thisrow{allreduce}/(\minvalue*\numiterations))] {\datafile};
	\legend{PM,mTTV,NLS,Factor Comm,Gram Comm};
\end{axis}
\end{tikzpicture}
\caption{Time breakdown for $243\times243\times243\times243$ tensor and rank $R=10$ on 81 processors for varying processor grids.}
\label{fig:commsweep10}
\end{figure}

In order to illustrate the effect of processor grid choice on running time, we show in \Cref{fig:commsweep10} a time breakdown over various processor grid choices for a 4D problem on 81 processors.
Because the tensor is cubical and 81 has a restricted factorization into 4 numbers, there are 5 distinct processor grids.
The overall takeaway is that the processor grid has very little effect on running time; in this experiment there is less than 10\% variation in overall time.
While the optimal processor grid reduced the communication time by approximately $3\times$ compared to the other processor grids, the running time is dominated by local computation, so it had little effect on overall time.
Furthermore, adjusting the processor grid affects the local tensor dimensions and the performance of the local computations, and the optimal processor grid led to slower local performance.
For $R=10$, all of the local computation is memory bandwidth bound, and we believe the variations in running times to be effects of some temporary quantities fitting into smaller levels of cache.

\subsection{Varying Approximation Rank}

One of the challenges of the CP (and NNCP) decomposition in practice is the choice of decomposition rank.
The most common technique is to compute multiple CP decompositions for various ranks.
As the rank $R$ increases, the approximation error  $\|\TA - \T{M}\|$ decreases with the better approximation power of more parameters. 
However, the benefit of increasing $R$ eventually diminishes if the data can be well approximated with a CP model.
Towards this end, we experiment with various values of $R$ to observe the relative increase in running time for two real-world data sets. 

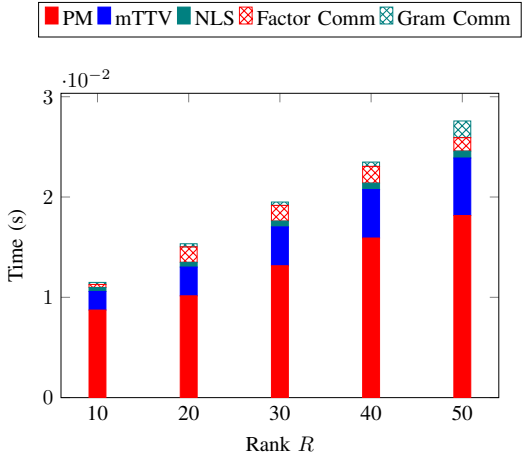
\begin{figure}
\renewcommand{\datafile}{ksw_4D_HSI.dat}
\renewcommand{\numiterations}{10}
\begin{tikzpicture}[scale=.8]
\begin{axis}[	
	ybar stacked,
	bar width=8pt,
	width=\columnwidth,
	height =.75\columnwidth,
	ylabel={Time (s)}, 
	xlabel={Rank $R$},
	y label style={yshift=-.5cm},
	ymin=0,
	symbolic x coords={D10,D20,D30,D40,D50},
	xticklabels={,10,20,30,40,50},
	legend style={at={(0.5,1.3)},anchor=north},
	legend columns=-1,
]
	\setcolors
	\addplot table[x=alg-K, y expr=(\thisrow{mttkrp}/(\minvalue*\numiterations))] {\datafile};
	\addplot table[x=alg-K, y expr=(\thisrow{mttv}/(\minvalue*\numiterations))] {\datafile};
	\addplot table[x=alg-K, y expr=((\thisrow{nnls}+\thisrow{gram})/(\minvalue*\numiterations))] {\datafile};
	\addplot table[x=alg-K, y expr=((\thisrow{allgather}+\thisrow{reducescatter})/(\minvalue*\numiterations))] {\datafile};
	\addplot table[x=alg-K, y expr=(\thisrow{allreduce}/(\minvalue*\numiterations))] {\datafile};
	\legend{PM,mTTV,NLS,Factor Comm,Gram Comm};
\end{axis}
\end{tikzpicture}
\caption{Per-iteration time breakdown of our implementation (using dimension trees) over various ranks for a time-lapse HSI dataset with dimensions $1344\times 1024\times 33 \times 9$ on 64 processors arranged in a $8\times8\times1\times1$ grid.}
\label{fig:ksweep4DHSI}
\end{figure}

\Cref{fig:ksweep4DHSI} shows the time breakdown of our implementation using a dimension tree on the 4D HSI dataset for $R=\{10,\dots,50\}$. 
We observe an overall time increase with increased $R$, but each part of the computation scales slightly differently.
The multi-TTV computation ({\em mTTV}) scales linearly with the increasing $R$,  whereas the partial MTTKRP ({\em PM}) is scaling super-linearly. 
This is because mTTV is cast as matrix-vector products (DGEMV) and PM is cast as matrix-matrix products (DGEMM).
As $R$ increases from 10 to 50, DGEMM performance improves but DGEMV performance is constant.  
The local NLS time is increasing with $O(R^3)$ as expected and the All-Reduce required for the Gram matrices scales with $O(R^2)$, becoming a significant cost for larger $R$. 

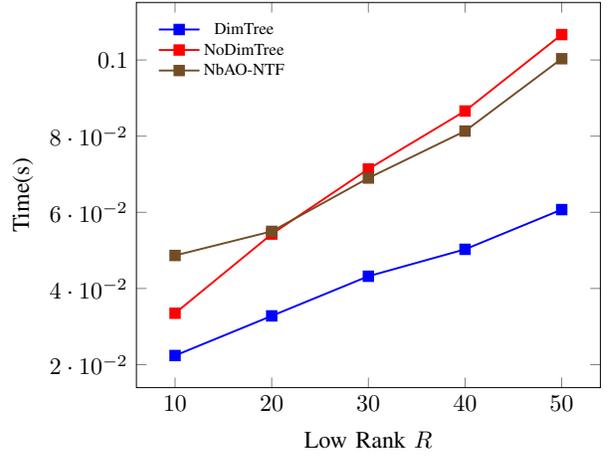
\begin{figure}
\begin{tikzpicture}[scale=.9]
\begin{axis}[legend style={draw=none, anchor=north,  cells={align=left}, nodes={scale=0.7}}, legend pos=north west, xlabel=Low Rank $R$, ylabel=Time(s), ylabel near ticks]
	\addplot+ [discard if not={alg}{D},thick,mark options={solid},mark=square*] table [x={k}, y expr=(\thisrow{total}/(\minvalue*\numiterations))] {kswp_neuro.dat};
	\addplot+ [discard if not={alg}{N},thick,mark options={solid},mark=square*] table [x={k}, y expr=(\thisrow{total}/(\minvalue*\numiterations))] {kswp_neuro.dat};
	\addplot+ [discard if not={alg}{L},thick,mark options={solid},mark=square*] table [x={k}, y expr=(\thisrow{total}/(\minvalue*\numiterations))] {kswp_neuro.dat};
	\legend{DimTree,NoDimTree,NbAO-NTF}
\end{axis}
\end{tikzpicture}
\caption{Overall running time for dFC dataset with dimensions $30{,}012\times 1200 \times 500$ on 1440 cores arranged in a $2\times 6\times 120$ processor grid and varying choices of rank $R$.}
\label{fig:ksweepneuro}
\end{figure}

In \Cref{fig:ksweepneuro} we compare performance for various ranks $R$ across all 3 algorithms, again used flat MPI. 
Starting at $R=10$ we see the largest speed up of $2\times$ for our implementation with a dimension tree over NbAO-NTF. 
This is due to a combination of the dimension tree performing fewer flops in the MTTKRPs and KRPs. 
However, as the rank increases this speed up diminishes to $1.6\times$. 
The loss of speed up is a result of the fact that, as we observed in \Cref{fig:ksweep4DHSI}, the multi-TTV operations do not scale as well as the partial-MTTKRPs for increasing $R$.
Again, the performance of our implementation without using dimension trees is comparable to NbAO-NTF. 

%% file: plots.tex

\newcommand{\datafile}{}
\newcommand{\alg}{}
\newcommand{\numiterations}{30}
\newcommand{\minvalue}{1}

\newif\ifnaive
\newif\ifksweep
\newif\iflegend
\newif\ifylabel
\newif\ifwider
\newif\ifliavas
\newif\ifweakthreeD

\newcommand{\setcolors}{
\pgfplotsset{cycle list={
	red, fill=red \\ 
	blue, fill=blue \\ 
	teal, fill=teal \\ 
	red, pattern=crosshatch, pattern color=red \\
	teal, pattern=crosshatch, pattern color=teal \\
	blue, pattern=crosshatch, pattern color=blue \\
}};
}

\newcommand{\breakdownplotoptions}{
	ybar stacked,
	reverse legend,
	bar width=8pt,
	ylabel={Time (s)}, 
	y label style={yshift=-.5cm},
	ymin=0,
	\ifweakthreeD
		symbolic x coords={D1,N1,D8,N8,D27,N27,D64,N64},
	\else
		symbolic x coords={D1,N1,D16,N16,D81,N81,D256,N256},
	\fi	
	xtick=data,
	legend style={at={(0.5,1.3)},anchor=north},
	legend columns=-1,
	reverse legend
}

\newcommand{\breakdownplot}{
\begin{axis}[\breakdownplotoptions]
	\setcolors
	\addplot table[x=alg-p, y expr=(\thisrow{mttkrp}/(\minvalue*\numiterations))] {\datafile};
	\addplot table[x=alg-p, y expr=(\thisrow{krp}/(\minvalue*\numiterations))] {\datafile};
	\addplot table[x=alg-p, y expr=((\thisrow{allgather}+\thisrow{reducescatter})/(\minvalue*\numiterations))] {\datafile};
	\addplot table[x=alg-p, y expr=((\thisrow{nnls}+\thisrow{gram}+\thisrow{allreduce})/(\minvalue*\numiterations))] {\datafile};
	\addplot table[x=alg-p, y expr=((\thisrow{err_compute}+\thisrow{err_communication})/(\minvalue*\numiterations))] {\datafile};
	\legend{MTTKRP,KRP,Factor Comm,NLS,Error};
\end{axis}
}

\newcommand{\labels}{
\node [align=center,text width=3cm] at (1cm, -1.15cm)   {\ifksweep 10 \else 1 \fi};
\node [align=center,text width=3cm] at (2.4cm, -1.15cm)   {\ifksweep 20 \else 16 \fi};
\node [align=center,text width=3cm] at (3.7cm, -1.15cm)   {\ifksweep 30 \else 81 \fi};
\node [align=center,text width=3cm] at (5.1cm, -1.15cm) {\ifksweep 40 \else 256 \fi};
}

\pgfplotsset{
    discard if not/.style 2 args={
        x filter/.append code={
            \edef\tempa{\thisrow{#1}}
            \edef\tempb{#2}
            \ifx\tempa\tempb
            \else
                \def\pgfmathresult{inf}
            \fi
        }
    }
}

\newcommand{\strongscalingplotoptions}{
		log basis y={2},
		log basis x={2},
		\ifliavas
			xlabel=Cores,
		\else
			xlabel=Nodes,
		\fi
		ylabel=Time (s),
		y tick label style={
	        /pgf/number format/.cd,
            	precision=4,
		},
	legend style={draw=none, cells={align=left}, nodes={scale=0.7}}
}

\newcommand{\strongscalingplot}{
\begin{loglogaxis}[\strongscalingplotoptions]
	\ifliavas
		\addplot+ [discard if not={alg}{DF},thick,mark options={solid},mark=square*] table [x={p}, y expr=(\thisrow{total}/(\minvalue*\numiterations))] {\datafile};
		\addplot+ [discard if not={alg}{NF},thick,mark options={solid},mark=square*] table [x={p}, y expr=(\thisrow{total}/(\minvalue*\numiterations))] {\datafile};
	\else
		\addplot+ [discard if not={alg}{D},thick,mark options={solid},mark=triangle*] table [x={p}, y expr=(\thisrow{total}/(\minvalue*\numiterations))] {\datafile};
		\addplot+ [discard if not={alg}{N},thick,mark options={solid},mark=square*] table [x={p}, y expr=(\thisrow{total}/(\minvalue*\numiterations))] {\datafile};
	\fi
	\ifliavas
		\addplot+ [discard if not={alg}{L},thick,mark options={solid}] table [x={p}, y expr=(\thisrow{total}/(\minvalue*\numiterations))] {\datafile};
		\legend{DimTree,NoDimTree,NbAO-NTF \cite{LK+17b}}
	\else
		\legend{DimTree,NoDimTree,FlatDimTree,FlatNoDimTree};
	\fi
\end{loglogaxis}
}

%% file: conclusion.tex
\section{Conclusion} \label{sec:conclusion}

In this work, we present a new implementation for distributed-memory NNCP that will be made publicly available.
The algorithm is general enough to handle any number of modes in the data tensor and can be adapted to use any NLS algorithm within the context of BCD (ALS).
We use a dimension tree optimization to avoid unnecessary recomputation within the bottleneck local MTTKRP computation, and we use an efficient parallelization that minimizes communication cost.
Our performance results show the ability to scale well to high processor counts, and we show favorable performance in comparison to state-of-the-art software for 3D tensors.

In particular, the performance results demonstrate that computing NNCP for dense tensors involves heavy computation relative to the sizes of the computed factor matrices.
By avoiding the communication of tensor entries and communicating only the factor matrices, the parallel algorithm is nearly always compute bound.
This observation is supported by the theoretical analysis: although the communication does not scale as well with $P$, the total amount of data depends on a sum of tensor dimensions rather the product of the tensor dimensions, which determines the total amount of computation.
For a relative comparison, consider memory-efficient parallel dense $n\times n$ matrix multiplication: the ratio of the $O(n^3/P)$ local computation to the $O(n^2/P^{1/2})$ communication is approximately the square root of the size of the local data, or $O((n^2/P)^{1/2})$.
In the case of NNCP, the ratio of local computation to communication is $O((I/P)^{1-1/N}/N)$, which is approximately the size of the local data raised to the power $1-1/N$.
This exponent is larger than $1/2$ and grows with $N$, and therefore it predicts the NNCP computation should be more computation bound than matrix multiplication.
Note that this analysis does not consider the type of local computation; for small $R$, the local computation will likely be memory bandwidth bound, but the algorithm will spend more time on local computation than on interprocessor communication.

We can also conclude from the performance results that the dimension tree optimization is the key to performance improvement over the state-of-the-art approaches.
For 3D tensors, we observe a benefit larger than the theoretical 50\% reduction in computation, and for larger numbers of modes, the improvement is only magnified.
Besides the reduction in flops, the dimension tree approach enjoys better DGEMM performance and avoids memory-bound KRP computations.
Furthermore, we see that tuning the processor grid had much less effect on overall performance.
Not only do reductions in communication not matter as much as computation, but different local tensor sizes can also cause variations in local performance that outweigh the savings in communication. 


%% file: appendix.tex

\appendix
\label{sec:appendix}

\begin{algorithm}
\caption{$(\CPl,\epsilon) = \text{Par-NNCP}(\TA,R)$}
\algfontsize
\label{alg:Par-NNCP-long}
\begin{algorithmic}[1]
\Require $\TA$ is an $I_1\times \cdots \times I_N$ tensor distributed across a $P_1\times \cdots \times P_N$ grid of $P$ processors, so that $\TA_{\V{p}}$ is $(I_1/P_1)\times \cdots \times (I_N/P_N)$ and is owned by processor $\V{p}=(p_1,\dots,p_N)$, $R$ is rank of approximation
\State \Comment{Initialize data}
\State $a = \text{Norm-Squared}(\TA_{\V{p}})$
\State $\alpha = \text{All-Reduce}(a,\textsc{All-Procs})$
\State $\epsilon = $ \texttt{Inf}
\For{$n=2$ to $N$}
	\State Initialize $\Mn{H}{n}_{\V{p}}$ of dimensions $(I_n/P)\times R$ 
	\State $\M[\overline]{G} = \text{Local-SYRK}(\Mn{H}{n}_{\V{p}})$
	\State $\Mn{G}{n} = \text{All-Reduce}(\M[\overline]{G},\textsc{All-Procs})$
	\State $\Mn{H}{n}_{p_n} = \text{All-Gather}(\Mn{H}{n}_{\V{p}},\textsc{Proc-Slice}(n,\VE{p}{n}))$
\EndFor
\State \Comment{Compute NNCP approximation}
\While{$\epsilon > $ \texttt{tol}}
	\State \Comment{Perform outer iteration of BCD}
	\For{$n=1$ to $N$}
	\State \Comment{Compute new factor matrix in $n$th mode}
	\State $\M[\overline]{M} = \text{Local-MTTKRP}(\TA_{p_1\cdots p_N},\{\Mn{H}{i}_{p_i}\},n)$
	\State $\Mn{M}{n}_{\V{p}} = \text{Reduce-Scatter}(\M[\overline]{M},\textsc{Proc-Slice}(n,\VE{p}{n}))$ 
	\State $\Mn{S}{n} = \Mn{G}{1} \Hada \cdots \Hada \Mn{G}{n-1} \Hada \Mn{G}{n+1} \Hada \cdots \Hada \Mn{G}{N}$
	\State $\Mn[\hat]{H}{n}_{\V{p}} = \text{Local-NLS-Update}(\Mn{S}{n},\Mn{M}{n}_{\V{p}})$
	\State \Comment{Normalize columns}
	\State $\V[\overline]{\lambda} = \text{Local-Col-Norms}(\Mn[\hat]{H}{n}_{\V{p}})$
	\State $\V{\lambda} = \text{All-Reduce}(\V[\overline]{\lambda},\textsc{All-Procs})$
	\State $\Mn{H}{n}_{\V{p}} = \text{Local-Col-Scale}(\Mn[\hat]{H}{n}_{\V{p}},\V{\lambda})$
	\State \Comment{Organize data for later modes}
	\State $\M[\overline]{G} = {\Mn{H}{n}_{\V{p}}}^\Tra\Mn{H}{n}_{\V{p}}$
	\State $\Mn{G}{n} = \text{All-Reduce}(\M[\overline]{G},\textsc{All-Procs})$
	\State $\Mn{H}{n}_{p_n} = \text{All-Gather}(\Mn{H}{n}_{\V{p}},\textsc{Proc-Slice}(n,\VE{p}{n}))$
	\EndFor
	\State \Comment{Compute relative error $\epsilon$ from mode-$N$ matrices}
	\State $\overline{\beta} = \text{Inner-Product}(\Mn{M}{N}_{\V{p}},\Mn[\hat]{H}{N}_{\V{p}})$
	\State $\beta = \text{All-Reduce}(\overline{\beta},\textsc{All-Procs})$
	\State $\gamma = \V{\lambda}^\Tra (\Mn{S}{N} \Hada \Mn{G}{N}) \V{\lambda}$
	\State $\epsilon = \sqrt{(\alpha-2\beta+\gamma)/\alpha}$ 
\EndWhile
\Ensure $\|\TA - \dsquare{\V{\lambda}; \Mn{H}{1},\dots,\Mn{H}{N}}\| /\|\TA\| = \epsilon$
\Ensure Local matrices: $\Mn{H}{n}_{\V{p}}$ is $(I_n/P)\times R$ and owned by processor $\V{p}=(p_1,\dots,p_N)$, for $1\leq n \leq N$, $\V{\lambda}$ stored redundantly on every processor
\end{algorithmic}
\end{algorithm}

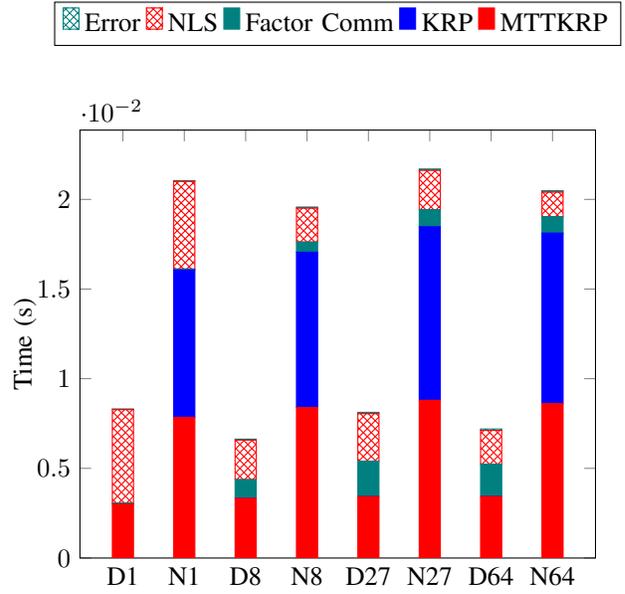
\begin{figure}
\begin{tikzpicture}
\renewcommand{\datafile}{wk_3D_syn2_scale_pgf.dat}
\renewcommand{\numiterations}{10}
\weakthreeDtrue
\breakdownplot
\end{tikzpicture}
\caption{Weak scaling of 3D synthetic tensors with (D) and without (N) the use of dimension trees.  The tensor and processor grid dimensions are $128k\times 128k\times 128k$ and $k\times k\times k$ for $k\in\{1,2,3,4\}$, and the rank is fixed at 32.  The reported times are per iteration.}
\label{fig:weaksynthetic3D}
\end{figure}